\documentclass[preprint,authoryear,12pt]{elsarticle}

\usepackage[english]{babel}
\usepackage[utf8]{inputenc}
\usepackage{fancyhdr}
\usepackage{amsthm}
\usepackage{amssymb}
\usepackage{amsmath}
\usepackage{enumerate}
\usepackage{algorithmic,algorithm}
\usepackage{url}
\usepackage{graphicx}
\usepackage{natbib}
\usepackage{lineno}

\newcommand{\EnsDisc}{\tau_d}
\newcommand{\Disc}[1]{\mathbf{#1}^{\EnsDisc}}
\newcommand{\R}[1]{\ensuremath{\mathbb{R}^{#1}}}
\newcommand{\learningSet}{S_{n,\EnsDisc}}
\newcommand{\Classifier}[1]{\phi_{n,\EnsDisc}^{#1}}
\newcommand{\proba}[1]{\mathbb{P}\left(#1\right)}
\newcommand{\Esp}[1]{\mathbb{E}\left(#1\right)}
\newcommand{\Sobolev}[1]{\mathcal{H}^m_{#1}}
\newcommand{\spline}[2]{\widehat{#1}_{#2,\EnsDisc}}
\newcommand{\Bound}[1]{\ensuremath{B^{#1}}}
\newcommand{\MatPS}[1]{\mathbf{M}_{#1,\EnsDisc}}
\newcommand{\inner}[3]{\langle #1,#2\rangle_{#3}}
\newcommand{\MatChol}[1]{\mathbf{Q}_{#1,\EnsDisc}}
\newcommand{\SplineOp}[1]{\mathcal{S}_{{#1},\EnsDisc}}
\newcommand{\InvSplineOp}[1]{\mathcal{S}^{-1}_{{#1},\EnsDisc}}
\newcommand{\Norm}[2]{\left\|#1\right\|_{#2}}
\newcommand{\Bayes}[1]{\ensuremath{L^*_{#1}}}

\newcommand{\BayesSplines}{\Bayes{d}}

\newtheorem{theorem}{Theorem}
\newtheorem{corollary}{Corollary}

\newtheorem{assump}{Assumption}

\theoremstyle{definition}
\newtheorem{remark}{Remark}

\journal{Pattern Recognition Letters}

\begin{document}
\sloppy

\begin{frontmatter}

\title{Consistency of Functional Learning Methods Based on Derivatives}

\author[fr]{Fabrice Rossi}
\ead{Fabrice.Rossi@telecom-paristech.fr}
\address[fr]{Télécom ParisTech, LTCI - UMR CNRS 5141, France}

\author[nvv1,nvv2]{Nathalie Villa-Vialaneix\corref{nvv-cor}}
\ead{nathalie.villa@math.univ-toulouse.fr}
\address[nvv1]{IUT de Perpignan (Dpt STID, Carcassonne), Université de Perpignan Via Domitia, France}
\address[nvv2]{Institut de Mathématiques de Toulouse, Université de Toulouse, France}
\cortext[nvv-cor]{Corresponding author.}

\begin{abstract}
In some real world applications, such as spectrometry, functional models achieve better predictive performances if they work on the derivatives of order $m$ of their inputs rather than on the original functions. As a consequence, the use of derivatives is a common practice in functional data analysis, despite a lack of theoretical guarantees on the asymptotically achievable performances of a derivative based model. In this paper, we show that a smoothing spline approach can be used to preprocess multivariate observations obtained by sampling functions on a discrete and finite sampling grid in a way that leads to a consistent scheme on the original infinite dimensional functional problem. This work extends \cite{mas_pumo_JNS2009} to nonparametric approaches and incomplete knowledge. To be more precise, the paper tackles two difficulties in a nonparametric framework: the information loss due to the use of the derivatives instead of the original functions and the information loss due to the fact that the functions are observed through a discrete sampling and are thus also unperfectly known: the use of a smoothing spline based approach solves these two problems. Finally, the proposed approach is tested on two real world datasets and the approach is experimentaly proven to be a good solution in the case of noisy functional predictors.
\end{abstract}

\end{frontmatter}
\begin{linenumbers}
\section{Introduction}

As the measurement techniques are developping, more and more data are high dimensional vectors generated by measuring a continuous process on a discrete sampling grid. Many examples of this type of data can be found in real world applications, in various fields such as spectrometry, voice recognition, time series analysis, etc.

Data of this type should not be handled in the same way as standard
multivariate observations but rather analysed as \emph{functional} data: each
observation is a function coming from an input space with infinite dimension,
sampled on a high resolution sampling grid. This leads to a large number of
variables, generally more than the number of observations. Moreover,
functional data are frequently smooth and generate highly correlated variables
as a consequence. Applied to the obtained high dimensional vectors, classical
statistical methods (e.g., linear regression, factor analysis) often lead to
ill-posed problems, especially when a covariance matrix has to be inverted
(this is the case, e.g., in linear regression, in discriminant analysis and
also in sliced inverse regression). Indeed, the number of observed values for
each function is generally larger than the number of functions itself and these values are often strongly correlated. As a consequence, when these data are considered as multidimensional vectors, the covariance matrix is ill-conditioned and leads to unstable and unaccurate solutions in models where its inverse is required. Thus, these methods cannot be directly
used. During past years, several methods have been adapted to that particular
context and grouped under the generic name of Functional Data Analysis (FDA)
methods. Seminal works focused on linear methods such as factorial analysis
(\cite{deville_AINSEE1974,dauxois_pousse_TE1976,besse_ramsay_P1986,james_hastie_sugar_B2000},
among others) and linear models
\cite{ramsay_dalzell_JRSS1991,cardot_etal_SPL1999,james_hastie_JRSSB2001}; a
comprehensive presentation of linear FDA methods is given in
\cite{ramsay_silverman_FDA1997,ramsay_silverman_AFDA2002}. More recently,
nonlinear functional models have been extensively developed and include
generalized linear models \cite{james_JRSSB2002,james_silverman_JASA2005}, kernel nonparametric regression \cite{ferraty_vieu_NPFDA2006}, Functional Inverse Regression \cite{ferre_yao_S2003}, neural networks \cite{rossi_conanguez_NN2005,rossi_etal_N2005}, $k$-nearest neighbors \cite{biau_etal_IEEETIT2005,laloe_SPL2008}, Support Vector Machines (SVM), \cite{rossi_villa_N2006}, among a very large variety of methods.

In previous works, numerous authors have shown that the derivatives of the
functions lead sometimes to better predictive performances than the functions themselves in inference tasks, as they provide information about the shape or the regularity of the function. In particular applications such as spectrometry \cite{ferraty_vieu_NPFDA2006,rossi_etal_N2005,rossi_villa_N2006}, micro-array data \cite{dejean_etal_EURASIPJBSB2007} and handwriting recognition \cite{williams_etal_ICANN2006,ba_burkhardt_IEEETPAMI2004}, these characteristics lead to accurate predictive models. But, on a theoretical point of the view, limited results about the effect of the use of the derivatives instead of the original functions are available: \cite{mas_pumo_JNS2009} studies this problem for a linear model built on the first derivatives of the functions. In the present paper, we also focus on the theoretical relevance of this common practice and extend \cite{mas_pumo_JNS2009} to nonparametric approaches and incomplete knowledge.

More precisely, we address the problem of the estimation of the conditional
expectation $\Esp{Y|X}$ of a random variable $Y$ given a functional random
variable $X$. $Y$ is assumed to be either real valued (leading to a regression
problem) or to take values in $\{-1,1\}$ (leading to a binary classification
problem). We target two theoretical difficulties. The first difficulty is the
potential information loss induced by using a derivative instead of the
original function: when one replaces $X$ by its order $m$ derivative
$X^{(m)}$, consistent estimators (such as kernel models
\cite{ferraty_vieu_NPFDA2006}) guarantee an asymptotic estimation of
$\Esp{Y|X^{(m)}}$ but cannot be used directly to address the original problem,
namely estimating $\Esp{Y|X}$. This is a simple consequence of the fact that
$X\mapsto X^{(m)}$ is not a one to one mapping. The second difficulty is
induced by sampling: in practice, functions are never observed exactly but
rather, as explained above, sampled on a discrete sampling grid. As a
consequence, one relies on approximate derivatives, $\widehat{X}_\tau^{(m)}$
(where $\tau$ denotes the sampling grid). This approach induces even more
information loss with respect to the underlying functional variable $X$: in
general, a consistent estimator of $\Esp{Y|\widehat{X}_\tau^{(m)}}$ will not
provide a consistent estimation of $\Esp{Y|X}$ and the optimal predictive
performances for $Y$ given $\widehat{X}_\tau^{(m)}$ will be lower than the
optimal predictive performances for $Y$ given $X$. 

We show in this paper that the use of a smoothing spline based approach solves both problems. Smoothing splines are used to estimate the functions from their sampled version in a convergent way. In addition, properties of splines are used to obtain estimates of the derivatives of the functions with no induced information loss. Both aspects are implemented as a preprocessing step applied to the multivariate observations generated via the sampling grid. The preprocessed observations can then be fed into any finite dimensional consistent regression estimator or classifier, leading to a consistent estimator for the original infinite dimensional problem (in real world applications, we instantiate the general scheme in the particular case of kernel machines \cite{shawetaylor_cristianini_KMPA2004}). 

The remainder of the paper is organized as follows: Section~\ref{notations} introduces the model, the main smoothness assumption and the notations. Section~\ref{splines} recalls important properties of spline smoothing. Section~\ref{consist} presents approximation results used to build a general consistent classifier or a general consistent regression estimator in Section~\ref{consist_general}. Finally, Section~\ref{application} illustrates the behavior of the proposed method for two real world spectrometric problems. The proofs are given at the end of the article.

\section{Setup and notations\label{notations}}

\subsection{Consistent classifiers and regression functions}

We consider a pair of random variables $(X,Y)$ where $X$ takes values in a functional space $\cal X$ and $Y$ is either a real valued random variable (regression case) or a random variable taking values in $\{-1,1\}$ (binary classification case). From this, we are given a learning set $S_n=\{(X_i,Y_i)\}_{i=1}^n$ of $n$ independent copies of $(X,Y)$. Moreover, the functions $X_i$ are not entirely known but sampled according to a non random sampling grid of finite length, $\EnsDisc=(t_l)_{l=1}^{|\EnsDisc|}$: we only observe $\Disc{X}_i=(X_i(t_1),\ldots X_i(t_{|\EnsDisc|}))^T$, a vector of $\R{|\EnsDisc|}$ and denote $\learningSet$ the corresponding learning set. Our goal is to construct:
\begin{enumerate}
	\item \emph{in the binary classification case}:  a classifier, $\Classifier{}$, whose misclassification probability
	\[
		L(\Classifier{})=\proba{\Classifier{}(\Disc{X})\neq Y}
	\]
	asymptotically reaches the Bayes risk
	\[
		\Bayes{}=\inf_{\phi:{\cal X}\rightarrow \{-1,1\}} \proba{\phi(X)\neq Y}
	\]
	i.e., $\lim_{|\EnsDisc|\rightarrow+\infty} \lim_{n\rightarrow+\infty} \Esp{L(\Classifier{})} = \Bayes{}$ ;
	\item \emph{in the regression case}: a regression function,
          $\Classifier{}$, whose $L^2$ error 
\[
L(\Classifier{})=\Esp{[\Classifier{}(\Disc{X})-Y]^2}
\]
asymptotically reaches the minimal $L^2$ error
	\[
		\Bayes{}=\inf_{\phi:{\cal X}\rightarrow \R{}} \Esp{[\phi{}(\Disc{X})-Y]^2}
	\]
i.e., $\lim_{|\EnsDisc|\rightarrow+\infty} \lim_{n\rightarrow+\infty}
L(\Classifier{}) = \Bayes{} $. 

This definition implicitly requires $\Esp{Y^2}<\infty$ and as a consequence,
corresponds to a $L^2$ convergence of $\Classifier{}$ to the conditional
expectation $\phi^*=\Esp{Y|X}$, i.e., to $\lim_{|\EnsDisc|\rightarrow+\infty}
\lim_{n\rightarrow+\infty} \Esp{[\Classifier{}(\Disc{X})-\phi^*(X)]^2}= 0$. 
\end{enumerate}
Such $\Classifier{}$ are said to be \emph{(weakly) consistent}
\cite{devroye_gyorfi_lugosi_PTPR1996,gyorfi_etal_DFTNR2002}. We have
deliberately used the same notations for the (optimal) predictive performances
in both the binary classification and the regression case. We will call
$\Bayes{}$ the Bayes risk even in the case of regression. Most of the
theoretical background of this paper is common to both the regression case and
the classification case: the distinction between both cases will be made only
when necessary. 

As pointed out in the introduction, the main difficulty is to show that the
performances of a model built on the $\Disc{X}_i$ asymptotically reach the
best performance achievable on the original functions $X_i$. In addition, we
will build the model on derivatives estimated from the $\Disc{X}_i$. 

\subsection{Smoothness assumption}
Our goal is to leverage the functional nature of the data by allowing differentiation operators to be applied to functions prior their submission to a more common classifier or regression function. Therefore we assume that the functional space $\cal X$ contains only differentiable functions. More precisely, $\cal X$ is the Sobolev space $\Sobolev{} =\Bigl\{ h\in L^2([0,1])\mid\,\forall\,j=1,\ldots,m,\ D^jh \text{ exists in the weak sense,}\,\text{and } D^mh\in L^2([0,1])\Bigr\}$, where $D^jh$ is the $j$-{\it th} derivative of $h$ (also denoted by $h^{(j)}$) and for an integer $m\geq 1$. Of course, by a straightforward generalization, any bounded interval can be considered instead of $[0,1]$.

To estimate the underlying functions $X_i$ and their derivatives from sampled data, we rely on smoothing splines. More precisely, let us consider a deterministic function $x \in \Sobolev{}$ sampled on the aforementioned grid. A smoothing spline estimate of $x$ is the solution, $\spline{x}{\lambda}$, of
\begin{equation}\label{eq:spline:definition}
	\arg\min_{h\in \Sobolev{}} \frac{1}{|\EnsDisc|} \sum_{l=1}^{|\EnsDisc|} (x(t_l)-h(t_l))^2 + \lambda \int_{[0,1]} (h^{(m)}(t))^2 dt,
\end{equation}
where $\lambda$ is a regularization parameter that balances interpolation
error and smoothness (measured by the $L^2$ norm of the $m$-{\it th}
derivative of the estimate). The goal is to show that a classifier or a
regression function built on $\spline{X}{\lambda}^{(m)}$ is consistent for the
original problem (i.e., the problem defined by the pair $(X,Y)$): this means
that using $\spline{X}{\lambda}^{(m)}$ instead of $X$ has no dramatic
consequences on the accuracy of the classifier or of the regression
function. In other words, asymptotically, no information loss occurs when one
replaces $X$ by $\spline{X}{\lambda}^{(m)}$. 

The proof is based on the following steps:
\begin{enumerate}
\item First, we show that building a classifier or a regression function on
  $\spline{X}{\lambda}^{(m)}$ is approximately equivalent to building a
  classifier or a regression function on
  $\Disc{X}=(X(t_l))_{l=1}^{|\EnsDisc|}$ using a specific metric. This is done
  by leveraging the Reproducing Kernel Hilbert Space (RKHS) structure of
  $\Sobolev{}$. This part serves one main purpose: it provides a solution to
  work 
  with estimation of the derivatives of the original function in a way that
  preserves all the information available in $\Disc{X}$. In other words, the
  best predictive performances for $Y$ theoretically available by building a
  multivariate model on $\Disc{X}$ are equal to the best predictive
  performances obtained by building a functional model on
  $\spline{X}{\lambda}^{(m)}$.  
\item Then, we link $\Esp{Y|\spline{X}{\lambda}}$ with $\Esp{Y|X}$ by
  approximation results available for smoothing splines. This part of the
  proof handles the effects of sampling. 
\item Finally, we glue both results via standard $\R{|\EnsDisc|}$ consistency
  results.
\end{enumerate}

\section{Smoothing splines and differentiation operators\label{splines}}

\subsection{RKHS and smoothing splines}

As we want to work on derivatives of functions from $\mathcal{H}^m$, a natural inner product for two functions of $\mathcal{H}^m$ would be $(u,v)\rightarrow \int_0^1 u^{(m)}(t) v^{(m)}(t) dt$. However, we prefer to use an inner product of $\mathcal{H}^m$ ($\int_0^1 u^{(m)}(t) v^{(m)}(t) dt$ only induces a semi-norm on $\mathcal{H}^m$) because, as will be shown later, such an inner product is related to an inner product between the sampled functions considered as vectors of \R{|\EnsDisc|}. 

This can be done by decomposing $\Sobolev{}$ into $\Sobolev{}=\Sobolev{0}\oplus\Sobolev{1}$ \cite{kimerldorf_wahba_JMAA1971}, where $\Sobolev{0}=\textrm{Ker} D^m=\mathbb{P}^{m-1}$ (the space of polynomial functions of degree less or equal to $m-1$) and $\Sobolev{1}$ is an infinite dimensional subspace of $\Sobolev{}$ defined via $m$ boundary conditions. The boundary conditions are given by a full rank linear operator from $\Sobolev{}$ to \R{m}, denoted $B$, such that $\textrm{Ker} B\cap \mathbb{P}^{m-1}=\{0\}$. Classical examples of boundary conditions include the case of ``natural splines'' (for $m=2$, $h(0)=h(1)=0$) and constraints that target only the first values of $h$ and its derivatives at a fixed position, for instance the conditions: $h(0)=\ldots=h^{(m-1)}(0)=0$. Other boundary conditions can be used \cite{berlinet_thomasagnan_RKHSPS2004,besse_ramsay_P1986,craven_wahba_NM1978}, depending on the application. 

Once the boundary conditions are fixed, an inner product on both $\Sobolev{0}$ and $\Sobolev{1}$ can be defined:
	\[
	\inner{u}{v}{1}=\inner{D^mu}{D^mv}{L^2}=\int_0^1 u^{(m)}(t) v^{(m)}(t) dt
	\]
is an inner product on $\Sobolev{1}$ (as $h\in \Sobolev{1}$ and $D^m h\equiv 0$ give $h\equiv 0$). Moreover, if we denote $B=(\Bound{j})_{j=1}^m$, then $\inner{u}{v}{0}=\sum_{j=1}^m \Bound{j}u \Bound{j}v$ is an inner product on $\Sobolev{0}$. We obtain this way an inner product on $\Sobolev{}$ given by
\begin{eqnarray*}
	\inner{u}{v}{\Sobolev{}}&=&\int_0^1 u^{(m)}(t)v^{(m)}(t) dt+\sum_{j=1}^m \Bound{j}u\Bound{j}v\\
	&=&\inner{\mathcal{P}^m_1(u)}{\mathcal{P}^m_1(v)}{1}+\inner{\mathcal{P}^m_0(u)}{\mathcal{P}^m_0(v)}{0}
\end{eqnarray*}
where $\mathcal{P}^m_i$ is the projector on $\Sobolev{i}$.

Equipped with $\inner{.}{.}{\Sobolev{}}$, $\Sobolev{}$ is a Reproducing Kernel Hilbert Space (RKHS, see e.g. \cite{berlinet_thomasagnan_RKHSPS2004,heckman_ramsay_CJS2000,wahba_SMOD1990}). More precisely, it exists a kernel $k: [0,1]^2\rightarrow \R{}$ such that, for all $u\in \Sobolev{}$ and all $t\in [0,1]$, $\inner{u}{k(t,.)}{\Sobolev{}}=u(t)$. The same occurs for $\Sobolev{0}$ and $\Sobolev{1}$ which respectively have reproducing kernels denoted by $k_0$ and $k_1$. We have $k=k_0+k_1$.

In the most common cases, $k_0$ and $k_1$ have already been explicitly calculated (see e.g., \cite{berlinet_thomasagnan_RKHSPS2004}, especially chapter 6, sections 1.1 and 1.6.2). For example, for $m\geq 1$ and the boundary conditions $h(0)=h'(0)=\ldots=h^{(m-1)}(0)=0$, we have: 
\[
	k_0(s,t)=\sum_{k=0}^{m-1} \frac{t^k s^k}{(k!)^2}.
\]
and
\[
	k_1(s,t)=\int_0^1 \frac{(t-w)_+^{m-1} (s-w)^{m-1}_+}{(m-1)!^2}\,dw.
\]

\subsection{Computing the splines}

We need now to compute to $\spline{x}{\lambda}$ starting with
$\Disc{x}=(x(t))_{t\in\EnsDisc}^T$. This can be done via a theorem from \cite{kimerldorf_wahba_JMAA1971}. We need the following compatibility assumptions between the sampling grid $\EnsDisc$ and the boundary conditions operator~$B$:
\begin{assump}\label{A_sampling_boundary}
	The sampling grid $\EnsDisc=(t_l)_{l=1}^{|\EnsDisc|}$ is such that
	\begin{enumerate}
		\item sampling points are distinct in $[0,1]$ and $|\EnsDisc|\geq m-1$
		\item the $m$ boundary conditions $\Bound{j}$ are linearly independent from the $|\EnsDisc|$ linear forms $h\mapsto h(t_l)$, for $l=1,\ldots,|\EnsDisc|$ (defined on $\Sobolev{}$)
	\end{enumerate}
\end{assump}
Then $\spline{x}{\lambda}$ and $\Disc{x}=(x(t))_{t\in\EnsDisc}^T$ are linked by the following result:
\begin{theorem}[\cite{kimerldorf_wahba_JMAA1971}]
	\label{th_kimeldorf_wahba}
	Under Assumption (A\ref{A_sampling_boundary}), the unique solution $\spline{x}{\lambda}$ to equation~\eqref{eq:spline:definition} is given by:
\begin{equation}\label{eq:spline:result}
\spline{x}{\lambda}=\SplineOp{\lambda}\Disc{x},
\end{equation}
where $\SplineOp{\lambda}$ is a full rank linear operator from $\R{|\EnsDisc|}$ to $\Sobolev{}$ defined by:
\begin{equation}
\label{eq:spline:operator}
\SplineOp{\lambda}=\omega^T M_0 + \eta^T M_1
\end{equation}
with
	\begin{itemize}
		\item $M_0=\left(U (K_{1}+\lambda I_d)^{-1} U^T\right)^{-1} U (K_1+\lambda I_d)^{-1}$
		\item $M_1=(K_{1}+\lambda I_d)^{-1}\left(I_d - U^TM_0\right)$;
		\item $\{\omega_1,\ldots,\omega_m\}$ is a basis of $\mathbb{P}^{m-1}$, $\omega=\left(\omega_1,\ldots,\omega_m\right)^T$ and $U=\left(\omega_i(t)\right)_{i=1,\ldots,m\ t\in\EnsDisc}$;
		\item $\eta=(k_1(t,.))_{t\in\EnsDisc}^T$  and $K_1=(k_1(t,t'))_{t, t'\in \EnsDisc}$.
	\end{itemize}
\end{theorem}

\subsection{No information loss}\label{subsectionNoLoss}
The first important consequence of Theorem \ref{th_kimeldorf_wahba} is that
building a model on $\spline{X}{\lambda}$ or on $\Disc{X}$ leads to the same
optimal predictive performances (to the same Bayes risk). This is formalized
by the following corollary:
\begin{corollary}\label{corollaryNoLoss}
Under Assumption (A\ref{A_sampling_boundary}), we have 
\begin{itemize}
\item in the binary classification case:
  \begin{equation}
    \label{eq:equiv_bayes}
		\begin{split}
   \inf_{\phi:\Sobolev{}\rightarrow \{-1,1\}}
   &\proba{\phi(\spline{X}{\lambda})\neq Y}=\\
   &\inf_{\phi:\R{|\EnsDisc|}\rightarrow \{-1,1\}} \proba{\phi(\Disc{X})\neq Y}
		\end{split}
 \end{equation}
\item in the regression case:
  \begin{equation}
    \label{eq:equiv_bayes_reg}
		\begin{split}
\inf_{\phi:\Sobolev{}\rightarrow \R{}} &\Esp{\left[\phi\left(\spline{X}{\lambda}\right)- Y\right]^2}=\\
&   \qquad \inf_{\phi:\R{|\EnsDisc|}\rightarrow \R{}} \Esp{\left[\phi\left(\Disc{X}\right)- Y\right]^2}
	\end{split}
  \end{equation}
\end{itemize}
\end{corollary}

\subsection{Differentiation operator\label{diff_norm}}
The second important consequence of Theorem \ref{th_kimeldorf_wahba} is that
the inner product $\inner{.}{.}{\Sobolev{}}$ is equivalent to a specific inner
product on \R{|\EnsDisc|} given in the following corollary: 
\begin{corollary}
	\label{cor_ps_h}
	Under Assumption (A\ref{A_sampling_boundary}) and for any $\Disc{u}=(u(t))_{t\in \EnsDisc}^T$ and $\Disc{v}=(v(t))_{t\in \EnsDisc}^T$ in \R{|\EnsDisc|}, 
	\begin{equation}
		\label{eq_ps_h}
		\inner{\spline{u}{\lambda}}{\spline{v}{\lambda}}{\Sobolev{}}=(\Disc{u})^T \MatPS{\lambda} \Disc{v}
	\end{equation}
	where $\MatPS{\lambda}= M_0^T W M_0 + M_1^T K_1 M_1$ with $W=(\inner{w_i}{w_j}{0})_{i,j=1,\ldots,m}$. The matrix $\MatPS{\lambda}$ is symmetric and positive definite and defines an inner product on \R{|\EnsDisc|}.
\end{corollary}
The corollary is a direct consequence of equations \eqref{eq:spline:result} and \eqref{eq:spline:operator}.

In practice, the corollary means that the euclidean space
$\left(\R{|\EnsDisc|},\inner{.}{.}{\MatPS{\lambda}}\right)$ is isomorphic to
$\left(\mathcal{I}_{\lambda,\EnsDisc},\inner{.}{.}{\Sobolev{}}\right)$,
where $\mathcal{I}_{\lambda,\EnsDisc}$ is the image of \R{|\EnsDisc|} by  
$\SplineOp{\lambda}$. As a consequence, one can use the Hilbert structure of
$\Sobolev{}$ directly in \R{|\EnsDisc|} via $\MatPS{\lambda}$: as the inner
product of $\Sobolev{}$ is defined on the order $m$ derivatives of the
functions, this corresponds to using those derivatives instead of the original
functions. 

More precisely, let $\MatChol{\lambda}$ be the transpose of the Cholesky
triangle of $\MatPS{\lambda}$ (given by the Cholesky decomposition
$\MatChol{\lambda}^T\MatChol{\lambda}=\MatPS{\lambda}$). Corollary
\ref{cor_ps_h} shows that $\MatChol{\lambda}$ acts as an approximate
differentiation operation on sampled functions.

Let us indeed consider an estimation method for multivariate inputs based only
on inner products or norms (that are directly derived from the inner products), such as, e.g., Kernel Ridge Regression
\cite{saunders_etal_ICML1998,shawetaylor_cristianini_KMPA2004}. In this latter
case, if a Gaussian kernel is used, the regression function has the following form:
\begin{equation}
	\label{eq::krrsolution}
u\mapsto\sum_{i=1}^nT_i\alpha_ie^{-\gamma \Norm{U_i-u}{\R{p}}^2}
\end{equation}
where $(U_i,T_i)_{1\leq i\leq n}$ are learning examples in
$\R{p}\times\{-1,1\}$ and the $\alpha_i$ are non negative real values obtained
by solving a quadratic programming problem and $\gamma$ is a parameter of the method. Then, if we use Kernel Ridge Regression on
the training set $\{(\MatChol{\lambda}\Disc{X}_i,Y_i)\}_{i=1}^n$ (rather than
the original training set $\{(\Disc{X}_i,Y_i)\}_{i=1}^n$), it will work on the
norm in $L^2$ of the derivatives of order $m$ of the spline estimates
of the $X_i$ (up to the boundary conditions). More precisely, the regression function will have
the following form:
\[
	\begin{split}
		\Disc{x}&\mapsto\sum_{i=1}^nY_i\alpha_ie^{-\gamma\Norm{\MatChol{\lambda}\Disc{X}_i-\MatChol{\lambda}\Disc{x}}{\R{|\EnsDisc|}}^2}\\
		&\mapsto\sum_{i=1}^nY_i\alpha_i e^{-\gamma\Norm{D^m
      \spline{X_i}{\lambda}-D^m \spline{x}{\lambda}}{L^2}^2}\\
		&\qquad \times e^{-\gamma\sum_{j=1}^m \left(\Bound{j}\spline{X_i}{\lambda}-\Bound{j}\spline{x}{\lambda}\right)^2}
	\end{split}
\]
In other words, up to the boundary
conditions, an estimation method based solely on inner products, or on norms
derived from these inner products, can be given
modified inputs that will make it work on an estimation of the derivatives of
the observed functions.  

\begin{remark}
As shown in Corollary \ref{corollaryNoLoss} in the previous section, building
a model on $\Disc{X}$ or on $\spline{X}{\lambda}$ leads to the same optimal
predictive performances. In addition, it is obvious that given any one-to-one
mapping $f$ from $\R{|\EnsDisc|}$ to itself, building a model on $f(\Disc{X})$
gives also the same optimal performances than building a model on
$\Disc{X}$. Then as $\MatChol{\lambda}$ is invertible, the optimal predictive
performances achievable with $\MatChol{\lambda}\Disc{X}$ are equal to the
optimal performances achievable with $\Disc{X}$ or with $\spline{X}{\lambda}$. 

In practice however, the actual preprocessing of the data can have a strong
influence on the obtained performances, as will be illustrated in Section
\ref{application}. The goal of the theoretical analysis of the present section
is to guarantee that no systematic loss can be observed as a consequence of
the proposed functional preprocessing scheme. 
\end{remark}

\section{Approximation results\label{consist}}
The previous section showed that working on $\Disc{X}$,
$\MatChol{\lambda}\Disc{X}$ or $\spline{X}{\lambda}$ makes no difference in
terms of optimal predictive performances. The present section addresses the
effects of sampling: asymptotically, the optimal predictive performances
obtained on $\spline{X}{\lambda}$ converge to the optimal performances
achievable on the original and unobserved functional variable $X$. 

\subsection{Spline approximation}

From the sampled random function $\Disc{X}=(X(t_1),\ldots,X(t_{|\EnsDisc|}))$, we can build an estimate, $\spline{X}{\lambda}$, of $X$. To ensure consistency, we must guarantee that $\spline{X}{\lambda}$ converges to $X$. In the case of a deterministic function $x$, this problem has been studied in numerous papers, such as 
\cite{craven_wahba_NM1978,ragozin_JAT1983,cox_SIAMJNA1984,utreras_JAT1988,wahba_SMOD1990} (among others). Here we recall one of the results which is particularly well adapted to our context. 

Obviously, the sampling grid must behave correctly, whereas the information contained in $\Disc{X}$ will not be sufficient to recover $X$. We need also the regularization parameter $\lambda$ to depend on $\EnsDisc$. Following \cite{ragozin_JAT1983}, a sampling grid $\EnsDisc$ is characterized by two quantities:
\begin{equation}
  \label{eq:Delta:Disc}
	\begin{split}
	&\overline{\Delta}_{\EnsDisc}=\max\{t_1,t_2-t_1,\ldots,1-t_{|\EnsDisc|}\}\\
	&\underline{\Delta}_{\EnsDisc}=\min_{1\leq i<|\EnsDisc|}\{t_{i+1}-t_i\}.
	\end{split}
\end{equation}
One way to control the distance between $X$ and $\spline{X}{\lambda}$ is to bound the ratio $\overline{\Delta}_{\EnsDisc}/\underline{\Delta}_{\EnsDisc}$ so as to ensure quasi-uniformity of the sampling grid. 

More precisely, we will use the following assumption:
\begin{assump}\label{A_spline_consistent}
There is R such that $\overline{\Delta}_{\EnsDisc}/\underline{\Delta}_{\EnsDisc}\leq R$ for all $d$. 
\end{assump}
Then we have:
\begin{theorem}[\cite{ragozin_JAT1983}]
	\label{th_cox}
	Under Assumptions (A\ref{A_sampling_boundary}) and (A\ref{A_spline_consistent}), there are two constants $A_{R,m}$ and $B_{R,m}$ depending only on $R$ and $m$, such that for any $x\in \Sobolev{}$ and any positive $\lambda$:
	\[
	\Norm{\spline{x}{\lambda}-x}{L^2}^2\leq \left(A_{R,m}\lambda+B_{R,m}\frac{1}{|\EnsDisc|^{2m}}\right)\Norm{D^mx}{L^2}^2.
	\]
\end{theorem}
This result is a rephrasing of Corollary 4.16 from \cite{ragozin_JAT1983} which is itself a direct consequence of Theorem 4.10 from the same paper.

Convergence of $\spline{x}{\lambda}$ to $x$ is then obtained by the following simple assumptions:
\begin{assump}\label{A_spline_consistent_more}
	The series of sampling points $\EnsDisc$ and the series of regularization parameters, $\lambda$, depending on $\EnsDisc$ and denoted by $(\lambda_d)_{d\geq 1}$, are such that $\lim_{d\rightarrow +\infty}|\EnsDisc|=+\infty$ and $\lim_{d\rightarrow+\infty}\lambda_d=0$.
\end{assump}

\subsection{Conditional expectation approximation}

The next step consists in relating the optimal predictive performances for the
regression and the classification problem $(X,Y)$ to the performances
associated to $(\spline{X}{\lambda_d},Y)$ when $d$ goes to infinity, i.e.,
relating $\Bayes{}$ to  
\begin{enumerate}
	\item \emph{binary classification case}:
		\[
		\BayesSplines=\inf_{\phi:\Sobolev{}\rightarrow \{-1,1\}} \proba{\phi(\spline{X}{\lambda_d})\neq Y},
		\]
	\item \emph{regression case}:
		\[
		\BayesSplines=\inf_{\phi:\Sobolev{}\rightarrow \R{}} \Esp{[\phi(\spline{X}{\lambda_d})- Y]^2}
		\]
\end{enumerate}

Two sets of assumptions will be investigated to provide the convergence of the Bayes risk $\BayesSplines$ to $\Bayes{}$:
\begin{assump}\label{A_esp_conditionnelle} {\bf Either}
	\begin{enumerate}[{(A\ref{A_esp_conditionnelle}}a{)}]
		\item \label{A_norme_X} $\Esp{\Norm{D^mX}{L^2}^2}$ is finite and $Y\in\{-1,1\}$,\\
		{\bf or} 
		\item \label{A_suite_disc} $\EnsDisc\subset \tau_{d+1}$ and $\Esp{Y^2}$ is finite.
	\end{enumerate}
\end{assump}
The first assumption (A\ref{A_esp_conditionnelle}\ref{A_norme_X}) requires an additional smoothing property for the predictor functional variable $X$ and is only valid for a binary classification problem whereas the second assumption (A\ref{A_esp_conditionnelle}\ref{A_norme_X}) requires an additional property for the sampling point series: they have to be growing sets.

Theorem~\ref{th_cox} then leads to the following corollary:
\begin{corollary}
	\label{cor_consist_d}
	Under Assumptions (A\ref{A_sampling_boundary})-(A\ref{A_esp_conditionnelle}), we have:
		\[
		\lim_{d\rightarrow +\infty} \BayesSplines = \Bayes{}.
		\]
\end{corollary}

\section{General consistent functional classifiers and regression functions\label{consist_general}}

\subsection{Definition of classifiers and regression functions on derivatives}\label{subsection:algorithm}

Let us now consider any consistent classification or regression scheme for standard multivariate data based either on the inner product or on the Euclidean distance between observations. Examples of such classifiers are Support Vector Machine \cite{steinwart_JC2002}, the kernel classification rule \cite{devroye_krzyzak_JSPI1989} and $k$-nearest neighbors \cite{devroye_gyorfi_NDELV1985,zhao_JMA1987} to name a few. In the same way, multilayer perceptrons \cite{lugosi_zeger_IEEETIT1995}, kernel estimates \cite{devroye_krzyzak_JSPI1989} and $k$-nearest neighbors regression \cite{devroye_etal_AS1994} are consistent regression estimators. Additional examples of consistent estimators in classification and regression can be found in \cite{devroye_gyorfi_lugosi_PTPR1996,gyorfi_etal_DFTNR2002}.

We denote $\psi_{\mathcal{D}}$ the estimator constructed by the chosen scheme using a dataset $\mathcal{D}=\{(U_i,T_i)_{1\leq i\leq n}\}$,  where the $(U_i,T_i)_{1\leq i\leq n}$ are $n$ independent copies of a pair of random variables $(U,T)$ with values in $\R{p}\times \{-1,1\}$ (classification) or $\R{p}\times \R{}$ (regression). 

The proposed functional scheme consists in choosing the estimator $\phi_{n,\EnsDisc}$ as $\psi_{\mathcal{E}_{n,\EnsDisc}}$ with the dataset $\mathcal{E}_{n,\EnsDisc}$ defined by:
\[
\mathcal{E}_{n,\EnsDisc}=\{(\MatChol{\lambda_d}\Disc{X}_i,Y_i)_{1\leq i\leq n}\}
\]
As pointed out in Section \ref{diff_norm}, the linear transformation
$\MatChol{\lambda_d}$ is an approximate multivariate differentiation operator:
up to the boundary conditions, an estimator based on
$\MatChol{\lambda_d}\Disc{X}$ is working on the $m$-th derivative of
$\spline{X}{\lambda_d}$.  

In more algorithmic terms, the estimator is obtained as follows:
\begin{enumerate}
\item choose an appropriate value for $\lambda_d$ 
\item compute $\MatPS{\lambda_d}$ using Theorem \ref{th_kimeldorf_wahba} and
  Corollary \ref{cor_ps_h};
\item compute the Cholesky decomposition of $\MatPS{\lambda_d}$ and the
  transpose of the Cholesky triangle, $\MatChol{\lambda_d}$ (such that
  $\MatChol{\lambda_d}^T\MatChol{\lambda_d}=\MatPS{\lambda_d}$);
\item compute  $\MatChol{\lambda_d}\Disc{X}_i$ to obtain the transformed
  dataset $\mathcal{E}_{n,\EnsDisc}$;
\item build a classifier/regression function $\psi_{\mathcal{E}_{n,\EnsDisc}}$
  with a multivariate method in $\R{|\EnsDisc|}$ applied to the dataset
  $\mathcal{E}_{n,\EnsDisc}$;
\item associate to a new sampled function $\Disc{X}_{n+1}$ the prediction
  $\psi_{\mathcal{E}_{n,\EnsDisc}}(\MatChol{\lambda}\Disc{X}_{n+1})$.
\end{enumerate}

\begin{figure}[htbp]
	\centering
	\includegraphics[width=0.95\linewidth]{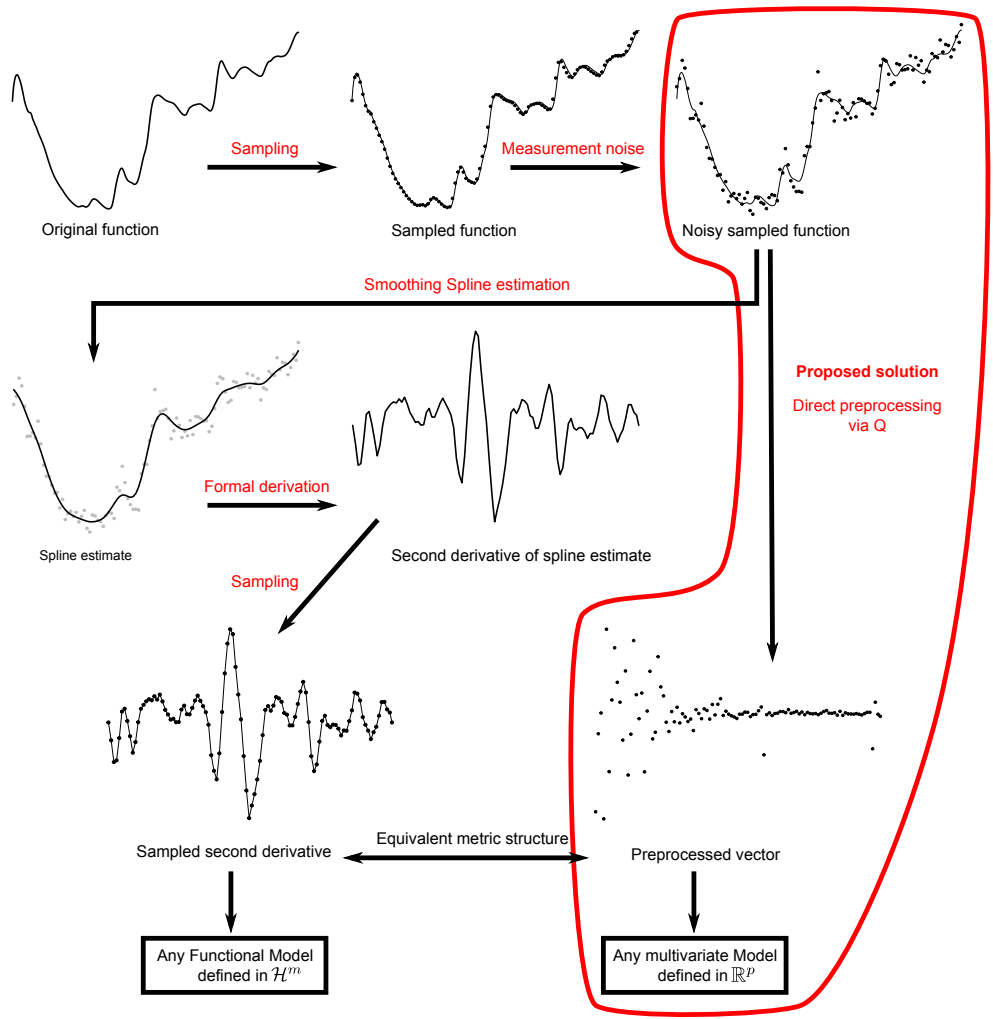}
	\label{fig:processing}
	\caption{Method scheme and its equivalence to the usual approach for using derivatives in learning algorithms.}
\end{figure}

Figure~\ref{fig:processing} illustrates the way the method performs: instead of relying on an approximation of the function and then on the derivation preprocessing of this estimates, it directly uses an equivalent metric by applying the $\MatChol{\lambda_d}$ matrix to the sampled function. The consistency result proved in Theorem~\ref{consistance_directe} shows that, combined with any consistent multidimensional learning algorithm, this method is (asymptotically) equivalent to using the original function drawn at the top left side of Figure~\ref{fig:processing}.

On a practical point of view, \cite{wahba_SMOD1990} demonstrates that cross
validated estimates of $\lambda$ achieve suitable convergence rates. Hence,
steps 1 and 2 can be computed simultaneously by minimizing the total cross
validated error for all the observations, given by
\[
	\sum_{i=1}^n \frac{1}{|\EnsDisc|} \sum_{t\in\EnsDisc} \frac{\left(x_i(t) - \spline{x_i}{\lambda}(t)\right)^2}{\left(1-A_{tt}(\lambda)\right)^2},
\]
where $A$ is a $|\EnsDisc|\times |\EnsDisc|$ matrix called the \emph{influence matrix} (see \cite{wahba_SMOD1990}), over a finite number of $\lambda$ values.

\subsection{Consistency result}
Corollary \ref{corollaryNoLoss} and Corollary~\ref{cor_consist_d} guarantee
that the estimator proposed in the previous section is consistent:  
\begin{theorem}
	\label{consistance_directe}
	Under assumptions (A\ref{A_sampling_boundary})-(A\ref{A_esp_conditionnelle}), the series of classifiers/regression functions $(\Classifier{})_{n,d}$ is consistent:
	\[
		\lim_{d\rightarrow+\infty}\lim_{n\rightarrow+\infty} \Esp{L\Classifier{}}=\Bayes{}
	\]
\end{theorem}

\subsection{Discussion}\label{discussion}
While Theorem \ref{consistance_directe} is very general, it could be easily extended to cover special cases such as additional hypothesis needed by the estimation scheme or to provide data based parameter selections. We discuss briefly those issues in the present section. 

It should first be noted that most estimation schemes, $\psi_{\mathcal{D}}$, depend on parameters that should fulfill some assumptions for the scheme to be consistent. For instance, in the Kernel Ridge Regression method in \R{p}, with Gaussian kernel, $\psi_{\mathcal{D}}$ has the form given in Equation~(\ref{eq::krrsolution}) where the $(\alpha_i)$ are the solutions of 
\[
	\begin{split}
	\arg\min_{\alpha\in\R{n}} &\sum_{i=1}^n \left(T_i - \sum_{j=1}^nT_j\alpha_je^{-\gamma\Norm{U_i-U_j}{\R{p}}^2}\right)^2 +\\ &\delta_n \sum_{i,j=1}^n T_iT_j \alpha_i\alpha_j e^{-\gamma\Norm{U_i-U_j}{\R{p}}^2}.
	\end{split}
\]
The method thus depends on the parameter of the Gaussian kernel, $\gamma$ and of the regularization parameter $\delta_n$. This method is known to be consistent if (see Theorem~9.1 of \cite{steinwart_christmann_SVM2008}):
\[
	\delta_n\xrightarrow{n\rightarrow+\infty}0\qquad \textrm{ and }\qquad n\delta_n^4 \xrightarrow{n\rightarrow+\infty} +\infty.
\]
Additional conditions of this form can obviously be directly integrated in Theorem~\ref{consistance_directe} to obtain consistency results specific to the corresponding algorithms. 

Moreover, practitioners generally rely on data based selection of the parameters of the estimation scheme $\psi_{\mathcal{D}}$ via a validation method: for instance, rather than setting $\delta_n$ to e.g., $n^{-5}$ for $n$ observations (a choice which is compatible with theoretical constraints on $\delta_n$), one chooses the value of $\delta_n$ that optimizes an estimation of the performances of the regression function obtained on an independent data set (or via a re-sampling approach). 

In addition to the parameters of the estimation scheme, functional data raise the question of the convenient order of the derivative, $m$, and of the sampling grid optimality. In practical applications, the number of available sampling points can be unnecessarily large (see \cite{biau_etal_IEEETIT2005} for an example with more than 8~000 sampling points). The preprocessing performed by $\MatChol{\lambda_d}$ do not change the dimensionality of the data which means that overfitting can be observed in practice when the number of sampling points is large compared to the number of functions. Moreover, processing very high dimensional vectors is time consuming. It is there quite interesting in practice to use a down-sampled version of the original grid. 

To select the parameters of $\psi_D$, the order of the derivative and/or the
down-sampled grid, a validation strategy, based on splitting the dataset into
training and validation sets, could be used. A simple adaptation of the idea
of
\cite{berlinet_etal_AISUP2008,biau_etal_IEEETIT2005,laloe_SPL2008,rossi_villa_N2006}
shows that a penalized validation method can be used to choose any combination
of those parameters consistently. According to those papers, the condition for the consistency of the validation strategy would simply relate the shatter coefficients of the set of classifiers in \R{d} to the penalization parameter of the validation. Once again, this type of results is a rather direct extension of Theorem \ref{consistance_directe}.

\section{Applications\label{application}}
In this section, we show that the proposed approach works as expected on real
world spectrometric examples: for some applications, the use of derivatives
leads to more accurate models than the direct processing of the spectra (see
e.g. \cite{rossi_etal_N2005,rossi_villa_N2006}
for other examples of such a behavior based on ad hoc estimators of the
spectra derivatives). It should be noted that the purpose of this section is
only to illustrate the behavior of the proposed method on finite datasets. The
theoretical results of the present paper show that all consistent schemes have
asymptotically identical performances, and therefore that using derivatives is
asymptotically useless. On a finite dataset however, preprocessing can have
strong influence on the predictive performances, as will be illustrated in
the present section. In addition, schemes that are not universally consistent,
e.g., linear models, can lead to excellent predictive performances on finite
datasets; such models are therefore included in the present section despite
the fact the theory does not apply to them. 

\subsection{Methodology}
The methodology followed for the two illustrative datasets is roughly the same:
\begin{enumerate}
\item the dataset is randomly split into a training set on which the
  model is estimated and a test set on which performances are
  computed. The split is repeated several times. The Tecator dataset (Section
  \ref{subsection:Tecator}) 
  is rather small (240 spectra) and exhibits a rather large
  variability in predictive performances between different random splits. We
  have therefore used 250 random splits. For the Yellow-berry dataset (Section \ref{subsection:yellowberry}), we used only 
  50 splits as the relative variability in performances is far less
  important. 
\item $\lambda$ is chosen by a global leave-one-out strategy on the
  spectra contained in training set (as suggested in Section
  \ref{subsection:algorithm}). More precisely, a leave-one-out
  estimate of the reconstruction error of the spline approximation of
  each training spectrum is computed for a finite set of candidate
  values for $\lambda$. Then a common $\lambda$ is chosen by
  minimizing the average over the training spectra of the leave-one-out
  reconstruction errors. This choice is relevant as cross validation
  estimates of $\lambda$ are known to have favorable theoretical
  properties (see \cite{craven_wahba_NM1978,utreras_SIAMJSC1981} among
  others). 
\item for regression problems, a Kernel Ridge Regression (KRR) \cite{saunders_etal_ICML1998,shawetaylor_cristianini_KMPA2004} is then performed to estimate the regression function; this method is consistent when used with a Gaussian kernel under additional conditions on the parameters (see Theorem~9.1 of \cite{steinwart_christmann_SVM2008}); as already explained, in the applications, Kernel Ridge Regression is performed both with a Gaussian kernel and with a linear kernel (in that last case, the model is essentially a
ridge regression model). Parameters of the models (a regularization parameter, $\delta_n$, in all cases and a
  kernel parameter, $\gamma$ for Gaussian kernels) are chosen by a grid search that
  minimizes a validation based estimate of the performances of the model (on
  the training set). A leave-one-out solution has
  been chosen: in Kernel Ridge Regression, the
  leave-one-out estimate of the performances of the model is obtained as a
  by-product of the estimation process, without additional computation cost,
  see e.g. \cite{cawley_talbot_NN2004}.\\
  Additionally, for a sake of comparison with a more traditional approach in
  FDA, Kernel Ridge Regression is compared with a nonparametric kernel
  estimate for the Tecator dataset
  (Section~\ref{subsec:tecator:regression}). Nonparametric kernel estimate is
  the first nonparametric approach introduced in Functional Data Analysis
  \cite{ferraty_vieu_NPFDA2006} and can thus be seen as a basis for comparison
  in the context of regression with functional predictors. For this method,
  the same methodology as with Kernel Ridge Regression was used: the parameter
  of the model (i.e., the bandwidth) was selected on a grid search minimizing
  a cross-validation estimate of the performances of the model. In this case,
  a 4-fold cross validation estimate was used instead of a leave-one-out
  estimate to avoid a large computational cost.

\item for the classification problem, a Support Vector Machine (SVM) is used
  \cite{shawetaylor_cristianini_KMPA2004}. As KRR, SVM are consistent when
  used with a Gaussian kernel \cite{steinwart_JC2002}. We also use a SVM with
  a linear kernel as this is quite adapted for classification in high
  dimensional spaces associated to sampled function data. We also use a
  K-nearest neighbor model (KNN) for reference. Parameters of the models (a
  regularization parameter for both SVM, a kernel parameter, $\gamma$ for
  Gaussian kernels and number of neighbors K for KNN) are chosen by a grid
  search that minimizes a validation based estimate of the classification
  error: we use a 4-fold cross-validation to get this estimate. 

\item We evaluate the models obtained for each random split on the test set. We report the mean and the standard deviation of the performance index (classification error and mean squared error, respectively) and assess the significance of differences between the reported figures via paired Student tests (with level 1\%). 

\item Finally, we compare models estimated on the raw spectra and on spectra transformed via the $\MatChol{\lambda_d}$ matrix for $m=1$ (first derivative) and $m=2$ (second derivative). For both values of $m$, we used the most classical boundary conditions ($x(0)=0$ and $Dx(0)=0$). Depending of the problem, other boundary conditions could be investigated but this is outside the scope of the present paper (see \cite{besse_ramsay_P1986,heckman_ramsay_CJS2000} for discussion on this subject). For the Tecator problem, we also compare these approaches with models estimated on first and second derivatives based on interpolating splines (i.e. with $\lambda=0$) and on first and second derivatives estimated by finite differences. \\
Note that the kind of preprocessing used has almost no impact on the computation time. In general, selecting the parameters of the model with
leave-one-out or cross-validation will use significantly more computing power
than constructing the splines and calculating their derivatives. For instance,
computing the optimal $\lambda$ with the approach described above
takes less than 0.1 second for the Tecator dataset on a standard PC using our
R implementation which is negligible compared to the several minutes used to
select the optimal parameters of the models used on the prepocessed data.
\end{enumerate}

\subsection{Tecator dataset}\label{subsection:Tecator}
The first studied dataset is the standard Tecator dataset \cite{thodberg_IEEETNN1995} \footnote{Data are
  available on statlib at \url{http://lib.stat.cmu.edu/datasets/tecator}}. It
consists in spectrometric data from the food industry. Each of the 240
observations is the near infrared absorbance spectrum of a meat sample
recorded on a Tecator Infratec Food and Feed Analyzer. Each spectrum is
sampled at 100 wavelengths uniformly spaced in the range 850--1050 nm. The
composition of each meat sample is determined by analytic chemistry and
percentages of moisture, fat and protein are associated this way to each
spectrum.

The Tecator dataset is a widely used benchmark in Functional Data Analysis, hence the motivation for its use for illustrative purposes. More precisely, in Section \ref{subsec:tecator:regression}, we address the original regression
problem by predicting the percentage of fat content from the
spectra with various regression method and various estimates of the derivative
preprocessing: this analysis shows that both the method and the use of
derivative have a strong effect on the performances whereas the way the
derivatives are estimated has almost no effect. Additionally, in Section
\ref{subsec:tecator:bruit}, we apply a noise (with various variances) to the
original spectra in order to study the influence of smoothing in the case of
noisy predictors: this section shows the relevance of the use of a smoothing
spline approach when the data are noisy. Finally,
Section~\ref{subsec:tecator:classif} deals with a classification problem
derived from the original Tecator problem (in the same way as what was done in
\cite{ferraty_vieu_CS2003}): conclusions of this section are similar to the
ones of the regression study.

\subsubsection{Fat content prediction}\label{subsec:tecator:regression}
As explained above, we first address the regression
problem that consists in predicting the fat content of peaces of meat from the Tecator dataset. The parameters of the model are optimized with a grid
search using the leave-one-out estimate of the predictive performances (both
models use a regularization parameter, with an additional width parameter in
the Gaussian kernel case). The original data set is split randomly into 160
spectra for learning and 80 spectra for testing. As shown in the result Table
\ref{table:TecatorReg}, the data exhibit a rather large variability; we use
therefore 250 random split to assess the differences between the different
approaches. 

The performance indexes are the mean squared error (M.S.E.) and the $R^2$.\footnote{$R^2=1-\frac{\textrm{M.S.E}}{\textrm{Var}(y)}$ where $\textrm{Var}(y)$ is the (empirical) variance of the target variable on the test set.} As a reference, the target variable (fat) has a variance equal to 14.36. Results are summarized in Table \ref{table:TecatorReg}.
\begin{table}[htbp]
  \begin{center}
	\begin{tabular}{llcc}
		Method & Data & Average M.S.E. & Average $R^2$\\
		&& and SD & \\
		\hline
		KRR Linear & O & 8.69 (4.47) & 95.7\%\\
		& S1 & 8.09 (3.85) & 96.1\%\\
		& IS1 & 8.09 (3.85) & 96.1\%\\
		& FD1 & 8.27 (4.17) & 96.0\%\\
		& S2 & 9.64 (4.98) & 95.3\%\\
		& IS2 & 9.87 (5.84) & 95.2\%\\
		& FD2 & 8.45 (4.18) & 95.9\%\\
		\hline
		KRR Gaussian & O & 5.02 (11.47) & 97.6\%\\
		& S1 & 0.485 (0.385) & 99.8\%\\
		& IS1 & 0.485 (0.385) & 99.8\%\\
		& FD1 & {\bf 0.484} (0.387) & {\bf 99.8\%}\\
		& S2 & 0.584 (0.303) & 99.7\%\\
		& IS2 & 0.586 (0.303) & 99.7\%\\
		& FD2 & 0.569 (0.281) & 99.7\%\\
		\hline
		NKE & O & 73.1 (16.5) & 64.2\%\\
		& S1 & 4.59 (1.09) & 97.7\%\\
		& IS1 & 4.59 (1.09) & 97.7\%\\
		& FD1 & 4.59 (1.09) & 97.7\%\\
		& S2 & 3.75 (1.22) & 98.2\%\\
		& IS2 & 3.75 (1.22) & 98.2\%\\
		& FD2 & 3.67 (1.18) & 98.2\%\\
		\hline
	\end{tabular}
  \caption{Summary of the performances of the chosen models on the test set
    (fat Tecator regression problem) when using either a kernel ridge regression (KRR) with linear kernel or with Gaussian kernel or when using a nonparametric kernel estimate (NKE) with various inputs: O (original data), S1 (smoothing splines with order 1 derivatives), IS1 (interpolating splines with order 1 derivatives), FD1 (order 1 derivatives estimated by finite differences) and S2, IS2 and FD2 (the same as previously with order 2 derivatives).}
\label{table:TecatorReg}
\end{center}
\end{table}

The first conclusion is that the method itself has a strong effect on the performances of the prediction: for this application, a linear method is not appropriate (mean squared errors are much greater for linear methods than for the kernel ridge regression used with a Gaussian kernel) and the nonparametric kernel estimate gives worse performances than the kernel ridge regression (indeed, they are about 10 times worse). Nevertheless, for nonparametric approaches (Gaussian KKR and NKE), the use of derivatives has also a strong impact on the performances: for kernel ridge regression, e.g., preprocessing by estimating the first order derivative leads to a strong decrease of the mean squared error. 

Differences between the average MSEs are not always significant, but we can
nevertheless rank the methods in increasing order of modeling error (using
notations explained in Table \ref{table:TecatorReg}) for Gaussian kernel ridge regression:
\[
\textrm{FD1} \leq \textrm{IS1} \leq \textrm{S1} < \textrm{DF2} \leq \textrm{SS2} < \textrm{IS2} < \textrm{O}
\]
where $<$ corresponds to a significant difference (for a paired Student test
with level 1\%) and $\leq$ to a non significant one. In this case, the data are very smooth and thus the use of smoothing splines instead of a finite differences approximation does not have a significant impact on the predictions. However, in this case, the roughest approach, consisting in the estimation of the derivatives by finite differences, gives the best performances.

\subsubsection{Noisy spectra}\label{subsec:tecator:bruit}
\begin{figure}[htbp]
  \centering
	\includegraphics[width=0.7\columnwidth,angle=270]{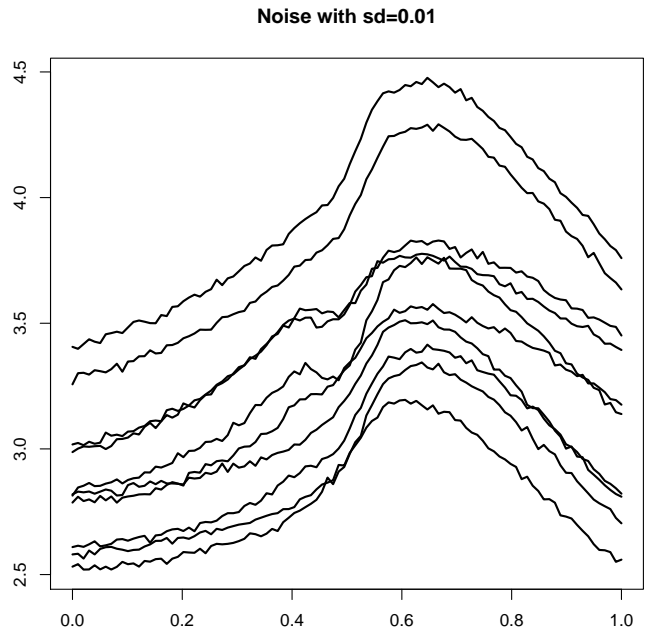}\\
	\includegraphics[width=0.7\columnwidth,angle=270]{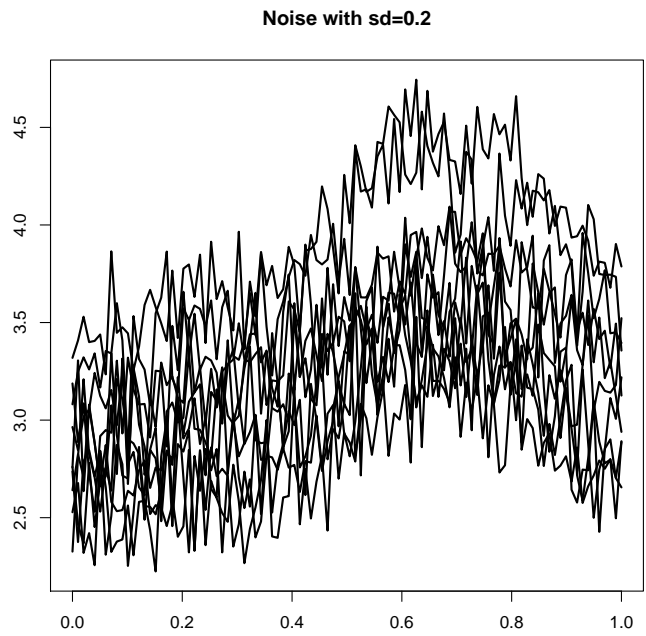}
  \caption{10 observations of the noisy data generated from the Tecator spectra as in Equation~\ref{eq::noise}}
  \label{fig::tecator-bruit}
\end{figure}

This section studies the situation in which functional data observations are
corrupted by noise. This is done by adding a noise to each spectrum of the
Tecator dataset. More precisely, each spectrum has been corrupted by 
\begin{equation}
	\label{eq::noise}
	X^b_i(t) = X_i(t) + \epsilon_{it} 
\end{equation}
where $(\epsilon_{it})$ are i.i.d. Gaussian variables with standard deviation equal to either 0.01 (small noise) or to 0.2 (large noise). 10 observations of the data generated this way are given in Figure~\ref{fig::tecator-bruit}.

\begin{figure}[htbp]
  \centering
	\includegraphics[width=0.7\columnwidth,angle=270]{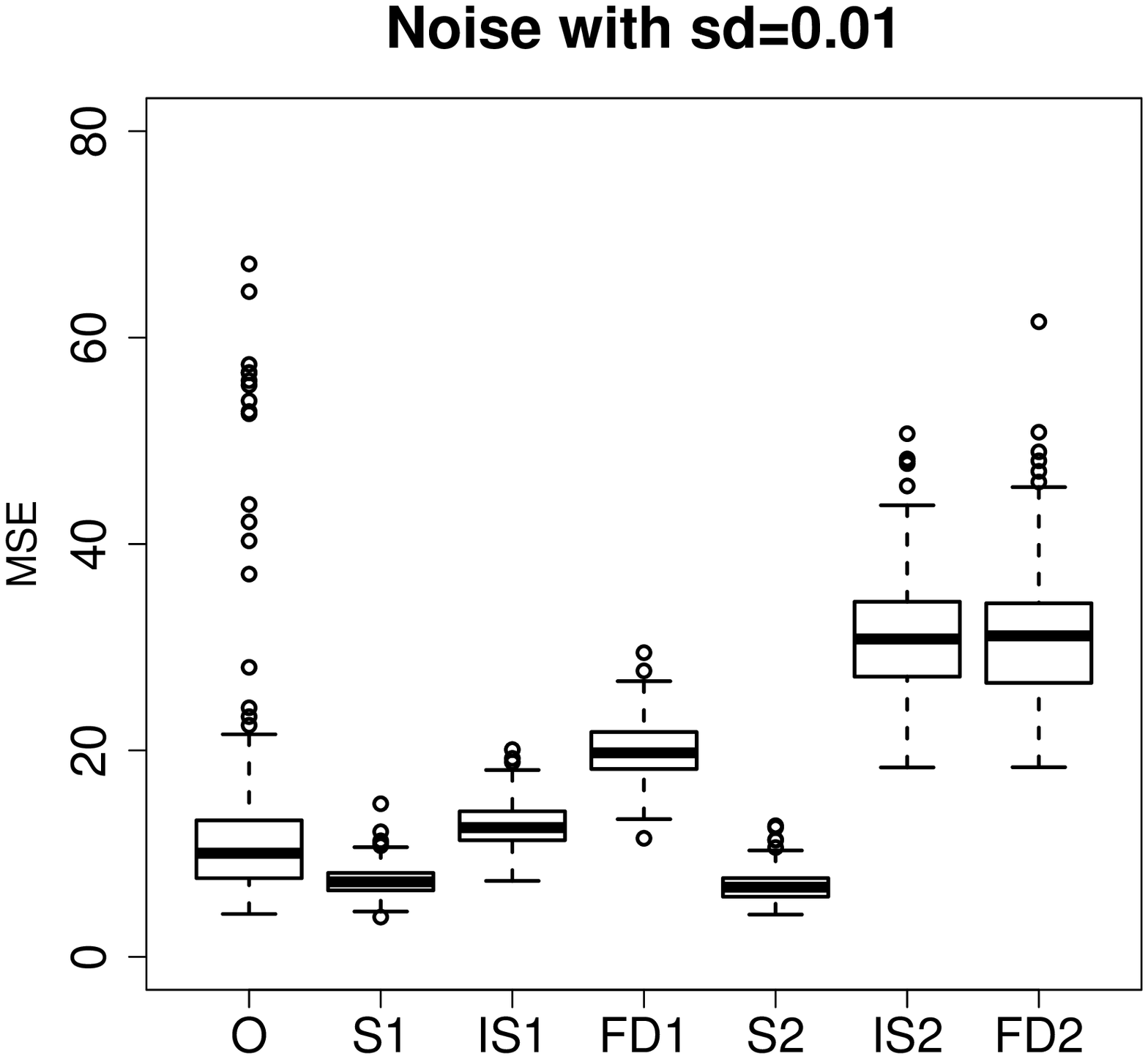}\\
	\includegraphics[width=0.7\columnwidth,angle=270]{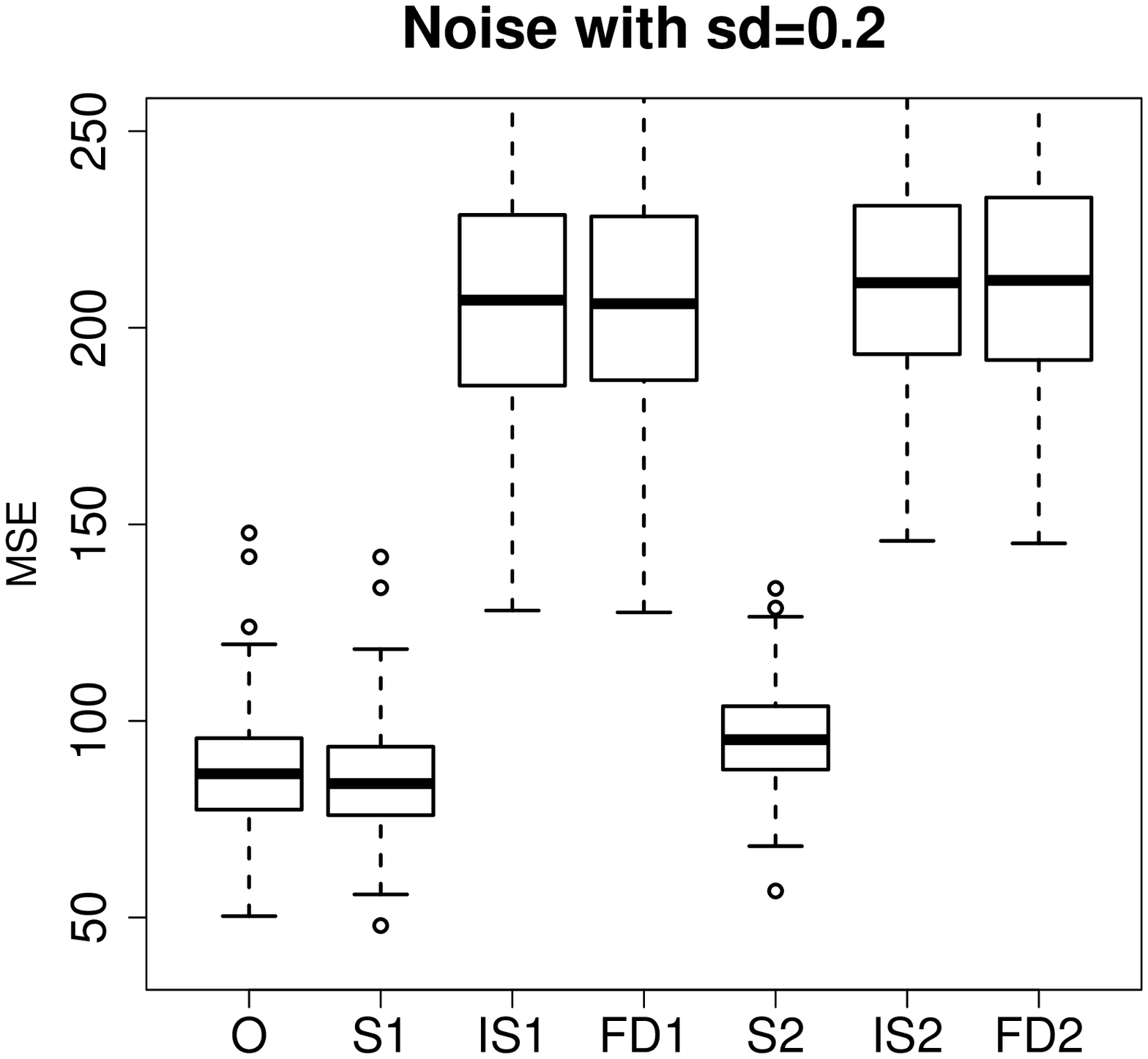}
  \caption{Mean squared errors boxplot for the noisy fat Tecator regression problem
    with Gaussian kernel (the worst test samples for IS and FD have been removed for a sake of clarity)}
  \label{fig:TecatorRegGaussian-bruit}
\end{figure}

The same methodology as for the non noisy data has been applied to $(X^b_i)$ to predict the fat content. The experiments have been restricted to the use of kernel ridge regression with a Gaussian kernel (according to the nonlinearity of the problem shown in the previous section). Results are summarized in Table \ref{table:TecatorReg-bruit} and Figure \ref{fig:TecatorRegGaussian-bruit}.
\begin{table}[htbp]
  \begin{center}
	\begin{tabular}{llcc}
		Noise & Data & Average M.S.E. & Average $R^2$\\
		& & and SD & \\
		\hline
		sd $=0.01$ & O & 13.3 (13.5) & 93.5\%\\
		& S1 & 7.45 (1.5) & 96.4\%\\
		& IS1 & 12.72 (2.2) & 93.8\%\\
		& FD1 & 20.03 (2.8) & 90.3\%\\
		& S2 & {\bf 6.83} (1.4) & {\bf 96.7\%}\\
		& IS2 & 31.23 (5.9) & 84.9\%\\
		& FD2 & 31.10 (5.9) & 84.9\%\\
		\hline
		sd $=0.2$ & O & 87.9 (13.9) & 57.4\%\\
		& S1 & {\bf 85.0} (12.5) & {\bf 58.8}\%\\
		& IS1 & 210.1 (36.1) & -1.9\%\\
		& FD1 & 209.1 (33.0) & -1.4\%\\
		& S2 & 95.9 (12.8) & 53.5\%\\
		& IS2 & 213.7 (33.1) & -3.6\%\\
		& FD2 & 235.1 (222.7) & -14.0\%\\
		\hline
	\end{tabular}
  \caption{Summary of the performances of the chosen models on the test set
    (fat Tecator regression problem) with noisy spectra.}
\label{table:TecatorReg-bruit}
\end{center}
\end{table}

In addition, the results can be ranked this way:
\begin{description}
	\item[Noise with sd equal to 0.01]
\[
	\textrm{S2} < \textrm{S1} < \textrm{IS1} \leq \textrm{O} < \textrm{FD1} < \textrm{IS2} \leq \textrm{FD2}
\]
	\item [Noise with sd equal to 0.2]
\[
	\textrm{S1} < \textrm{O} < \textrm{S2} < \textrm{FD1} < \textrm{IS1} < \textrm{IS2} \leq \textrm{FD2}
\]
\end{description}
where $<$ corresponds to a significant difference (for a paired Student test with level 1\%).

The first conclusion of these experiments is that, even though the derivatives
are the relevant predictors, their performances are strongly affected by the
noise (compared to the ones of the original data: note that the average
M.S.E. reported in Table \ref{table:TecatorReg} are more 10 times lower
that the best ones from Table \ref{table:TecatorReg-bruit} and that, in the best cases, $R^2$ is slightly greater than 50\% for the most noisy dataset). In particular, using interpolating splines or finite difference derivatives leads to highly deteriorated performances. In this situation, the approach proposed in the paper is particularly useful and helps to keep better performances than with the original data. Indeed, the differences of the smoothing splines approach with the original data is still significant (for both derivatives in the ``small noise'' case and for the first order derivative in the ``high noise'' case), even though, the most noisy the data are, the most difficult it is to estimate the derivatives in an accurate way. That is, except for smoothing spline derivatives, the estimation of the derivatives for the most noisy dataset is so bad that it leads to negative $R^2$ when used in the regression task.

\subsubsection{Fat content classification}\label{subsec:tecator:classif}
In this section, the fat content regression problem is transformed into a
classification problem. To avoid imbalance in class sizes, the median value of
the fat in the dataset is used as the splitting criterion: the first class
consists in 119 samples with strictly less than 13.5 \% of fat, while the second
class contains the other 121 samples with a fat content equal or higher than
13.5~\%. 

As in previous sections, the analysis is conducted on 250 random splits of the
dataset into 160 learning spectra and 80 test spectra. We used stratified
sampling: the test set contains 40 examples from each class. The
4 fold cross-validation used to select the parameters of the models on the
learning set is also stratified with roughly 20 examples of each class in each
fold. 

The performance index is the mis-classification rate (MCR) on the test set,
reported in percentage and averaged over the 250 random splits. Results are
summarized in Table \ref{table:TecatorClassif}.  As in the previous sections,
both the model and the preprocessing have some influence on the results. In
particular, using derivatives always improves the classification accuracy
while the actual method used to compute those derivatives has no particular
influence on the results. Additionally, using interpolation splines leads, in
this particular problem, to results that are exactly identical to the ones
obtained with the smoothing splines: they are not reported in Table
\ref{table:TecatorClassif}.

\begin{table}[htbp]
  \begin{center}
	\begin{tabular}{llcc}
		Method & Data & Average MCR & SD of MCR\\
		\hline
		Linear SVM& O &1.41  & 1.55\\
		& S1 & \textbf{0.73}& 1.15\\
		& FD1 & 0.74 & 1.15\\
		& S2 &  0.94& 1.27\\
		& FD2 & 0.92 & 1.23\\
		\hline
		Gaussian SVM& O & 3.39 & 2.57\\
		& S1 & 0.97& 1.41\\
		& FD1 & 0.98 & 1.42\\
		& S2 &  0.99& 2.00\\
		& FD2 & 0.97 & 1.27\\
		\hline
		KNN & O & 22.0 & 5.02\\
		& S1 & 6.67& 2.55\\
		& FD1 & 6.57 & 2.55\\
		& S2 &  1.93& 1.65\\
		& FD2 & 1.93 & 1.63\\
		\hline
	\end{tabular}
  \caption{Summary of the performances of the chosen models on the test set
    (Tecator fat classification problem). See Table \ref{table:TecatorReg} for
  notations. MCR stands for mis-classification rate, SD for standard deviation.}
\label{table:TecatorClassif}
\end{center}
\end{table}

More precisely, for the three models (linear SVM, Gaussian SVM and KNN),
differences in mis-classification rates between the smoothing spline
preprocessing and the finite differences calculation is never significant,
according to a Student test with level 1 \%. Additionally while the actual
average  mis-classification rates might seem quite different, the large
variability of the results (shown by the standard deviations) leads to
significant differences only for the most obvious cases. In particular, SVM
models using derivatives (of order one or two) are indistinguishable one from
another using a Student test with level 1 \%: all methods with less than 1 \%
of mean mis-classification rate perform essentially identically. Other
differences are significant: for instance the linear SVM used on raw data
performs significantly worse than any SVM model used on derivatives. 

It should be noted that the classification task studied in the present section
is obviously simpler than the regression task from which it is derived. This
explains the very good predictive performances obtained by simple models such
as a linear SVM, especially with the proper preprocessing. 

\subsection{Yellow-berry dataset}\label{subsection:yellowberry}
\begin{figure}[ht]
	\centering
	\includegraphics[angle=-90,width=0.95\columnwidth]{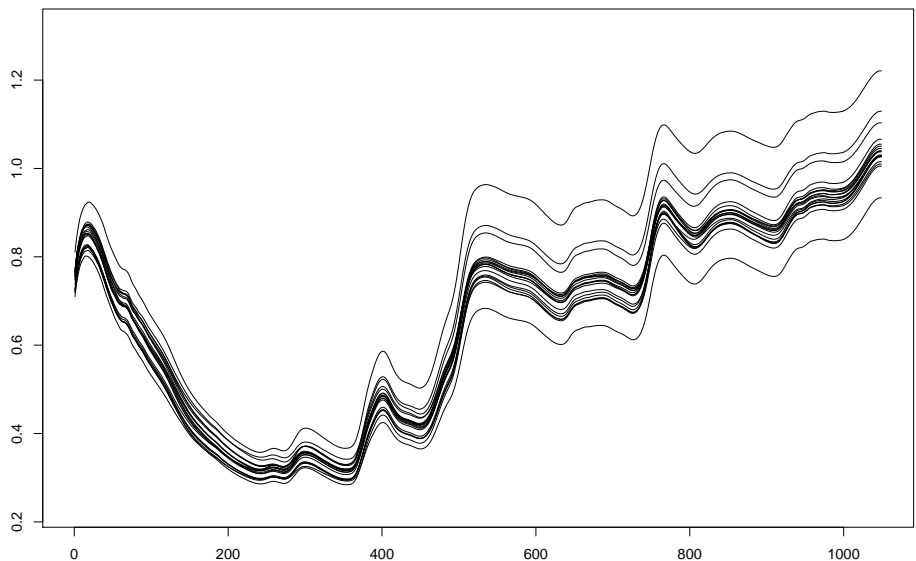}  
	\caption{20 observations of NIR spectra of durum wheat}
	\label{bledur}
\end{figure}

The goal of the last experiment is to predict the presence of yellow-berry in durum wheat (\emph{Triticum durum}) kernels via a near infrared spectral analysis (see Figure~\ref{bledur}). Yellow-berry is a defect of the durum wheat seeds that reduces the quality of the flour produced from affected wheat. The traditional way to assess the occurrence of yellow-berry is by visual analysis of a sample of the seed stock. In the current application, a quality measure related to the occurrence of yellow-berry is predicted from the spectrum of the seed. 

The dataset consists in 953 spectra sampled at 1049 wavelengths uniformly
spaced in the range 400--2498 nm. The dataset is split randomly into 600
learning spectra and 353 test spectra. Comparatively to the Tecator dataset,
the variability of the results is smaller in the present case. We used therefore 50 random splits rather than 250 in the previous section. 

The regression models were build via a Kernel Ridge Regression approach using
a linear kernel and a Gaussian kernel. In both cases, the regularization
parameter of the model is optimized by a leave-one-out approach. In addition,
the width parameter of the Gaussian kernel is optimized via the same procedure
at the same time.

\begin{figure}[htbp]
  \centering
\includegraphics[angle=270,width=0.95\columnwidth]{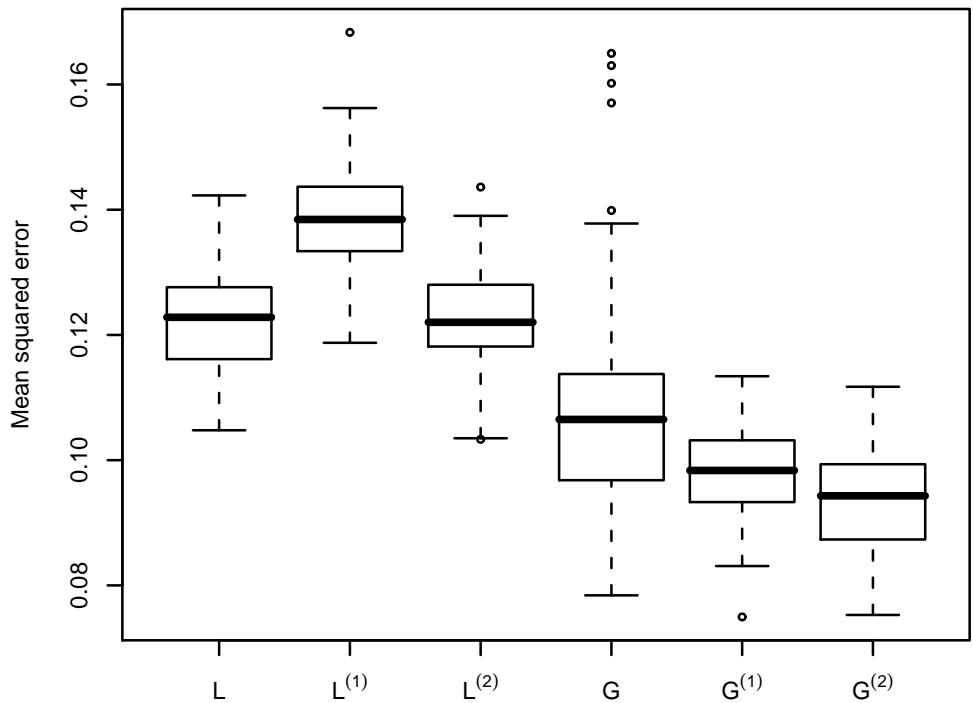}  
  \caption{Mean squared error boxplots for the ``durum wheat'' regression problem (see Table \ref{table:TriticumDurum} for the full names of the regression models)}
  \label{fig:TriticumDurum}
\end{figure}

The performance index is the mean squared error (M.S.E.). As a reference, the target variable has a variance of $0.508$. Results are summarized in Table \ref{table:TriticumDurum}  and Figure \ref{fig:TriticumDurum}. 
\begin{table}[htbp]
  \begin{center}
	\begin{tabular}{lccc}
		Kernel and Data & Average M.S.E. & Standard deviation&Average $R^2$\\\hline
		Linear-O & 0.122 & $8.77\, 10^{-3}$&76.1\%\\
		Linear-S1 & 0.138&$9.53\, 10^{-3}$&73.0\%\\
		Linear-S2 & 0.122&$8.41\, 10^{-3}$&76.1\%\\\hline
		Gaussian-O & 0.110&$20.2\, 10^{-3}$&78.5\%\\
		Gaussian-S1 &0.0978&$7.92\, 10^{-3}$&80.9\%\\
		Gaussian-S2 &0.0944&$8.35\, 10^{-3}$&81.5\%\\\hline
	\end{tabular}
  \caption{Summary of the performances of the chosen models on the test set (durum wheat regression problem)}
\label{table:TriticumDurum}
\end{center}
\end{table}

As in the previous section, we can rank the methods in increasing order of modelling error, we obtain the following result:
\[
\textrm{G-S2} < \textrm{G-S1} < \textrm{G-O} < \textrm{L-O} \leq \textrm{L-S2} < \textrm{L-S1},
\]
where G stands for Gaussian kernel and L for  linear kernel (hence G-S2 stands for kernel ridge regression with gaussian kernel and smoothing splines with order 2 derivatives); $<$ corresponds to a significant difference (for a paired Student test with level 1\%) and $\leq$ to a non significant one. For this application, there is a significant gain in using a non linear model (the Gaussian kernel). In addition, the use of derivatives leads to less contrasted performances that the ones obtained in the previous section but it still improves the quality of the non linear model in a significant way. In term of normalized mean squared error (mean squared error divided by the variance of the target variable), using a non linear model with the second derivatives of the spectra corresponds to an average gain of more than 5\% (i.e., a reduction of the normalised mean squared error from 24\% for the standard linear model to 18.6\%).

\section{Conclusion}

In this paper we proposed a theoretical analysis of a common practice that
consists in using derivatives in classification or regression problems when
the predictors are curves. Our method relies on smoothing splines
reconstruction of the functions which are known only via a discrete
deterministic sampling. The method is proved to be consistent for very general
classifiers or regression schemes: it reaches asymptotically the best risk
that could have been obtained by constructing a regression/classification
model on the true random functions. 

We have validated the approach by combining it with nonparametric regression
and classification algorithms to study two real-world spectrometric
datasets. The results obtained in these applications confirm once again that
relying on derivatives can improve the quality of predictive models compared
to a direct use of the sampled functions. The way the derivatives are
estimated does not have a strong impact on the performances except when the
data are noisy. In this case, the use of smoothing splines is quite relevant.

In the future, several issues could be addressed. An important practical
problem is the choice of the best order of the derivative, $m$. We consider
that a model selection approach relying on a penalized error loss could be
used, as is done, in e.g., \cite{rossi_villa_N2006}, to select the dimension
of truncated basis representation for functional data. Note that in practice,
such parameter selection method could lead to select $m=0$ and therefore to
automatically exclude derivative calculation when it is not needed. This will
extend the application range of the proposed model. 

A second important point to study it the convergence rate for the method. It
would be very convenient for instance, to be able to relate the size of the
sampling grid to the number of functions. But, this latter issue would require
the use of additional assumptions on the smoothness of the regression function
whereas the result presented in this paper, even if more limited, only needs
mild conditions.

\section{Acknowledgement}

We thank Cécile Levasseur and  Sylvain Coulomb (École d'Ingénieurs de Purpan, EIP, Toulouse, France) for sharing the interesting problem presented in Section~\ref{subsection:yellowberry}.

We also thank Philippe Besse (Institut de Mathématiques de Toulouse, Université de Toulouse, France) for helpfull discussions and suggestions.

Finally, we thank the anonymous reviewers for their valuable comments and suggestions that helped to improve the quality of the paper.

\section{Proofs}

\subsection{Theorem~\ref{th_kimeldorf_wahba}}

In the original theorem (Lemma~3.1) in \cite{kimerldorf_wahba_JMAA1971}, one has to verify that $(k_0(t_l,.))_l$ spans $\Sobolev{0}$ and that $(k_1(t_l,.))_l$ are linearly independent. These are consequences of Assumption (A\ref{A_sampling_boundary}).
	
	First, $k_0(s,t)=\sum_{i,j=0}^{m-1}   b_{ij}^{(-1)}s^it^j$ where $\widetilde{B}=(b_{i,j}^{(-1)})_{i,j}$ is the inverse of $(\sum_{l=1}^m\Bound{l}s^i\Bound{l}t^j)_{i,j}$ (see \cite{heckman_ramsay_CJS2000}). Then $(k_0(t_1,s),\ldots,k_0(t_{|\EnsDisc|},s))=(1,s,\ldots,s^{m-1})\widetilde{B}[V_{m-1}(t_1,\ldots,t_{|\EnsDisc|})]^T$ where $V_{m-1}(t_1,\ldots,t_{|\EnsDisc|})$ is the Vandermonde matrix with $m-1$ columns and $|\EnsDisc|$ rows associated to values $t_1,\ldots,t_{|\EnsDisc|}$. If the $(t_l)_l$ are distinct, this matrix is of full rank.

	Moreover the reproducing property shows that $\sum_{l=1}^{|\EnsDisc|} a_l k_1(t_l,.)\equiv 0$ implies $\sum_{l=1}^{|\EnsDisc|} a_l f(t_l)\equiv 0$ for all $f\in\Sobolev{1}$. Hence, $\Sobolev{1}=\textrm{Ker} \left(B^T,\sum_{l=1}^{\EnsDisc}a_l \zeta_l\right)^T$ where $\zeta_l$ denotes the linear form $h\in \Sobolev{}\rightarrow h(t_l)$. As the co-dimension of $\Sobolev{1}$ is $\textrm{dim} \Sobolev{0}= m$ and as, by Assumption (A\ref{A_sampling_boundary}), $B$ is linearly independent of $\sum_{l=1}^{\EnsDisc}a_l \zeta_l$, we thus have $\sum_{l=1}^{\EnsDisc}a_l \zeta_l\equiv 0$ (or $\textrm{codim}\,\textrm{Ker}\left(B^T,\sum_{l=1}^{\EnsDisc}a_l \zeta_l\right)^T = \textrm{dim}\, \textrm{Im} \left(B^T,\sum_{l=1}^{\EnsDisc}a_l \zeta_l\right)$ would be $m+1$). Thus, we obtain that $\sum_{l=1}^{|\EnsDisc|} a_l f(t_l)\equiv 0$ for all $f$ in $\Sobolev{}$ and, as $(t_l)$ are distinct, that $a_l=0$ for all $l$, leading to the independence conclusion for the  $(k_1(t_l,.))_l$.
	
	Finally, we prove that $\SplineOp{\lambda}$ is of full rank. Indeed, if $\SplineOp{\lambda}\Disc{x}=0$, $\omega^T M_0\Disc{x}=0$ and $\eta^T M_1\Disc{x}=0$. As $(\omega_k)_k$ is a basis of $\Sobolev{0}$, $\omega^T M_0\Disc{x}=0$ implies $M_0\Disc{x}=0$ and therefore $M_1=(K_{1}+\lambda I_d)^{-1}$.  As shown above, the $(k_1(t_l,.))_l$ are linearly independent and therefore $\eta M_1 \Disc{x}=0$ implies $M_1\Disc{x}=0$, which in turns leads to $\Disc{x}=0$ via the simplified formula for $M_1$.

\subsection{Corollary~\ref{corollaryNoLoss}}
We give only the proof for the classification case, the regression case is
identical. 

According to Theorem \ref{th_kimeldorf_wahba}, there is a full rank linear
mapping from \R{|\EnsDisc|} to $\Sobolev{}$, $\SplineOp{\lambda}$, such that
for any function $x\in\Sobolev{}$,
$\spline{x}{\lambda}=\SplineOp{\lambda}\Disc{x}$.  Let us denote
$\mathcal{I}_{\lambda,\EnsDisc}$ the image of \R{|\EnsDisc|} by
$\SplineOp{\lambda}$, $\mathbf{P}_{\lambda,\EnsDisc}$ the orthogonal
projection from $\Sobolev{}$ to $\mathcal{I}_{\lambda,\EnsDisc}$ and
$\InvSplineOp{\lambda}$ the inverse of $\SplineOp{\lambda}$ on
$\mathcal{I}_{\lambda,\EnsDisc}$. Obviously, we have $\InvSplineOp{\lambda}\circ \mathbf{P}_{\lambda,\EnsDisc}(\spline{x}{\lambda})=\Disc{x}$.

Let $\psi$ be a measurable function from \R{|\EnsDisc|} to $\{-1,1\}$. Then $\zeta_\psi$ defined on $\Sobolev{}$ by $\zeta_\psi(u)=\psi\left(\InvSplineOp{\lambda}\circ \mathbf{P}_{\lambda,\EnsDisc}(u)\right)$ is a measurable function from $\Sobolev{}$ to $\{-1,1\}$ (because $\InvSplineOp{\lambda}$ and $\mathbf{P}_{\lambda,\EnsDisc}$ are both continuous). Then for any measurable $\psi$, $\inf_{\phi:\Sobolev{}\rightarrow \{-1,1\}} \proba{\phi(\spline{X}{\lambda})\neq Y}\leq \proba{\zeta_\psi(\spline{X}{\lambda})\neq Y}=\proba{\psi(\Disc{X})\neq Y}$, 
and therefore
\begin{equation}\label{ineq_equiv_bayes_sobo}
	\begin{split}
 \inf_{\phi:\Sobolev{}\rightarrow \{-1,1\}} &\proba{\phi(\spline{X}{\lambda})\neq Y} \leq\\
 &\inf_{\phi:\R{|\EnsDisc|}\rightarrow \{-1,1\}} \proba{\phi(\Disc{X})\neq Y}. 
	\end{split}
\end{equation}
Conversely, let $\psi$ be a measurable function from $\Sobolev{}$ to $\{-1,1\}$. Then $\zeta_\psi$ defined on $\R{|\EnsDisc|}$ by $\zeta_\psi(\mathbf{u})=\psi(\SplineOp{\lambda}(\mathbf{u}))$, is measurable. Then for any measurable $\psi$, $\inf_{\phi:\R{|\EnsDisc|}\rightarrow \{-1,1\}} \proba{\phi(\Disc{X})\neq Y}\leq \proba{\zeta_\psi(\Disc{X})\neq Y}=\proba{\psi(\spline{X}{\lambda})\neq Y}$, and therefore
\begin{equation}\label{ineq_equiv_bayes_disc}
	\begin{split}\
		inf_{\phi:\R{|\EnsDisc|}\rightarrow \{-1,1\}}& \proba{\phi(\Disc{X})\neq Y}\leq\\
 &\inf_{\phi:\Sobolev{}\rightarrow \{-1,1\}} \proba{\phi(\spline{X}{\lambda})\neq Y}.
	\end{split}
\end{equation}
The combination of equations \eqref{ineq_equiv_bayes_sobo} and
\eqref{ineq_equiv_bayes_disc} gives equality \eqref{eq:equiv_bayes}.

\subsection{Corollary~\ref{cor_consist_d}}
\begin{enumerate}
		\item {\bf Suppose assumption (A\ref{A_esp_conditionnelle}\ref{A_norme_X}) is fullfilled}\\
		The proof is based on Theorem 1 in \cite{farago_gyorfi_IEEETIT1975}. This theorem relates the Bayes risk of a classification problem based on $(X,Y)$ with the Bayes risk of the problem $(T_d(X),Y)$ where $(T_d)$ is a series of transformations on $X$. 

		More formally, for a pair of random variables $(X,Y)$, where $X$ takes values in $\mathcal{X}$, an arbitrary metric space, and $Y$ in $\{-1,1\}$, let us denote for any series of functions $T_d$ from $\mathcal{X}$ to itself, $\Bayes{}(T_d)=\inf_{\phi:\,\mathcal{X}\rightarrow
			\{-1,1\}}\proba{\phi(T_d(X))\neq Y}$. Theorem 1 from \cite{farago_gyorfi_IEEETIT1975} states that $\Esp{\delta(T_d(X),X)}\xrightarrow{d\rightarrow +\infty} 0$ implies $\Bayes{}(T_d)\xrightarrow{d\rightarrow+\infty} \Bayes{}$, where $\delta$ denotes the metric on $\mathcal{X}$. 

		This can be applied to $\mathcal{X}=(\Sobolev{},\inner{.}{.}{L^2})$ with $T_d(X)=\spline{X}{\lambda_d}=S_{\lambda_d, \tau_d}\Disc{X}$: under Assumptions (A\ref{A_sampling_boundary}) and (A\ref{A_spline_consistent}), Theorem~\ref{th_cox} gives: $\Norm{T_d(X)-X}{L^2}^2\leq  \left(A_{R,m}\lambda_d+B_{R,m}\frac{1}{|\EnsDisc|^{2m}}\right)\Norm{D^mX}{L^2}^2$.
Taking the expectation of both sides gives $\Esp{\Norm{T_d(X)-X}{L^2}}\leq \left(A_{R,m}\lambda_d+B_{R,m}\frac{1}{|\EnsDisc|^{2m}}\right)\Esp{\Norm{D^mX}{L^2}^2}$, 
using the fact that the constants are independent of the function under analysis. Then under Assumptions (A\ref{A_esp_conditionnelle}\ref{A_norme_X}) and (A\ref{A_spline_consistent_more}), $\Esp{\Norm{T_d(X)-X}{L^2}}\xrightarrow{d\rightarrow+\infty}0$. According to \cite{farago_gyorfi_IEEETIT1975}, this implies $\lim_{d\rightarrow\infty}\BayesSplines=\Bayes{}.$ 
		
		\item {\bf Suppose assumption (A\ref{A_esp_conditionnelle}\ref{A_suite_disc}) is fullfilled}\\
		The conclusion will follow both for classification case and
                for regression case. The proof follows the general ideas of
                \cite{biau_etal_IEEETIT2005,rossi_conanguez_NPL2006,rossi_villa_N2006,laloe_SPL2008}. Under
                assumption~(A\ref{A_sampling_boundary}), by
                Theorem~\ref{th_kimeldorf_wahba} and with an argument similar to
                those developed in the proof of  Corollary \ref{corollaryNoLoss}, $\sigma(\spline{X}{\lambda_d})=\sigma(\{X(t)\}_{t\in\EnsDisc})$. From assumption (A\ref{A_esp_conditionnelle}\ref{A_suite_disc}), $\sigma(\{X(t)\}_{t\in\EnsDisc})$ is clearly a filtration. Moreover, as $\Esp{Y}$ and thus $\Esp{Y^2}$ are finite, $\Esp{Y|\spline{X}{\lambda_d}}$ is a uniformly bounded martingal for this filtration (see Lemma 35 of \cite{pollard_UGMTP2002}). This martingale converges in $L^1$-norm to $\Esp{Y|\sigma\left(\cup_d\sigma(\spline{X}{\lambda_d})\right)}$; we have
		\begin{itemize}
			\item $\sigma\left(\cup_d\sigma(\spline{X}{\lambda_d})\right) \subset \sigma(X)$ as $\spline{X}{\lambda_d}$ is a function of $X$ (via Theorem~\ref{th_kimeldorf_wahba});
			\item by Theorem~\ref{th_cox}, $\spline{X}{\lambda_d}\xrightarrow{d\rightarrow+\infty,\ \textrm{surely}}X$ in $L^2$ which proves that $X$ is $\sigma\left(\cup_d\sigma(\spline{X}{\lambda_d})\right)$-measurable.
		\end{itemize}
		Finally, $\Esp{Y|\sigma\left(\cup_d\sigma(\spline{X}{\lambda_d})\right)}=\Esp{Y|X}$ and $\Esp{Y|\spline{X}{\lambda_d}}\xrightarrow{d\rightarrow +\infty,\ L^1}\Esp{Y|X}$.
		
		The conclusion follows from the fact that:
		\begin{enumerate}
			\item {\it binary classification case:} the bound $\BayesSplines-\Bayes{} \leq 2\Esp{\left|\Esp{Y|\spline{X}{\lambda_d}}-\Esp{Y|X}\right|}$ (see Theorem~2.2 of \cite{devroye_gyorfi_lugosi_PTPR1996}) concludes the proof;
			\item {\it regression case:} as $\Esp{Y^2}$ is finite, $\Esp{\Esp{Y|\spline{X}{\lambda_d}}^2}$ is also finite and the convergence also happens for the quadratic norm (see Corollary 6.22 in \cite{kallenberg_FMP1997}), i.e.,
				\[
				\lim_{d\rightarrow+\infty}\Esp{\left(\Esp{Y|X}-\Esp{Y|\spline{X}{\lambda_d}}\right)^2} = 0
				\]
			Hence, as $\BayesSplines-\Bayes{}=\Esp{\left(\Esp{Y|X}-\Esp{Y|\spline{X}{\lambda_d}}\right)^2}$, the conclusion follows.
		\end{enumerate}
	\end{enumerate}

\subsection{Theorem~\ref{consistance_directe}}

We have
	\begin{equation}
		\label{dec_L}
		L(\phi_{n,d})-\Bayes{}=L\Classifier{}-\BayesSplines+\BayesSplines-\Bayes{}.
	\end{equation}
	Let $\epsilon$ be a positive real. By Corollary~\ref{cor_consist_d}, it exists $d_0 \in \mathbb{N}^*$ such that, for all $d\geq d_0$,
	\begin{equation}
		\label{util_cor_consist_d}
		\BayesSplines-\Bayes{} \leq \epsilon.
	\end{equation}
Moreover, as shown in Corollary \ref{corollaryNoLoss} and as
$\MatChol{\lambda_d}$ is invertible, we have in the binary
classification case: $\BayesSplines=\inf_{\phi:\R{|\EnsDisc|}\rightarrow \{-1,1\}}
\proba{\phi(\Disc{X})\neq Y}=\inf_{\phi:\R{|\EnsDisc|}\rightarrow \{-1,1\}}
\proba{\phi\left(\MatChol{\lambda_d}\Disc{X}\right)\neq Y}$, and in the regression case: $\BayesSplines=\inf_{\phi:\R{|\EnsDisc|}\rightarrow \R{}}
\Esp{\left[\phi\left(\Disc{X}\right)- Y\right]^2}=
\inf_{\phi:\R{|\EnsDisc|}\rightarrow \R{}}
\Esp{\left[\phi\left(\MatChol{\lambda_d}\Disc{X}\right)- Y\right]^2}$.
By hypothesis, for any fixed $d$, $\phi_{n,\tau_{d}}$ is consistent, that
is 
\[
\lim_{n\rightarrow+\infty}
\Esp{L(\phi_{n,\tau_{d}})}=\inf_{\phi:\R{|\EnsDisc|}\rightarrow \{-1,1\}} 
\proba{\phi\left(\MatChol{\lambda_{d}}\Disc{X}\right)\neq Y},
\]
in the classification case and
\[
\lim_{n\rightarrow+\infty}
\Esp{L(\phi_{n,\tau_{d}})}=\inf_{\phi:\R{|\EnsDisc|}\rightarrow \R{}}
\Esp{\left[\phi\left(\MatChol{\lambda_d}\Disc{X}\right)- Y\right]^2},
\]
in the regression case, and therefore for any fixed $d_0$,
$\lim_{n\rightarrow+\infty} \Esp{L(\phi_{n,\tau_{d_0}})}=L_{d_0}^*$. Combined
with equations \eqref{dec_L} and \eqref{util_cor_consist_d}, this concludes the
proof. 
\end{linenumbers}

\end{document}